\documentclass[fleqn]{article}

\newcommand{\titel} {A negative result on algebraic specifications of the meadow of rational numbers}

\usepackage{graphicx,version,stmaryrd,amssymb,amsmath,amsthm}
\usepackage[usenames,dvipsnames]{color}

\newtheorem{theorem}{Theorem}
  
\newtheorem{proposition}{Proposition}  
  
\newtheorem{corollary}{Corollary}

\theoremstyle{definition}
\newtheorem*{example}{Example}

\newcommand{\Md}{\mathsf{Md}}

\title{\titel}

\author{
	Jan A.\ Bergstra \& Inge Bethke 
}

\begin{document}
\maketitle

\begin{abstract}
$\mathbb{Q}_0$|the involutive meadow of the rational numbers|is the  {\mbox{} zero-totalized expansion field} of the rational numbers where the multiplicative inverse operation is made total by imposing $0^{-1}=0$.
In this note, we prove that $\mathbb{Q}_0$ cannot be specified  by the usual axioms for meadows augmented  by a finite set of axioms of the form $(1+  \cdots +1+x^2)\cdot (1+  \cdots +1 +x^2)^{-1}=1$. \\
\end{abstract}

\section{Introduction}
 {\mbox{} This note contributes to the algebraic datatype specification of $\mathbb{Q}_0$|the rational numbers equipped with a total inverse operation.
Advantages and disadvantages of dividing by zero in various ways have been amply discussed in the mathematics and logic literature (see e.g.\  \cite{K75,O83}) and we do not wish to add to those matters here.  The same holds for issues regarding the origins of various approaches to division by zero. But we believe that the observation made in \cite{BT2007} that $\mathbb{Q}_0$ has a finite initial algebra specification was at that time original. Here we elaborate on this theme.}

\noindent In \cite{BT2007} it is shown that 
\[
\mathbb{Q}_0 \cong \mathcal{I}(\Sigma_{\Md} , \Md+ L_4)
\]
where  $\mathcal{I}(\Sigma_{\Md} , \Md+ L_4)$ is the initial algebra of the theory $\Md$  of meadows given in Table \ref{meadow} enriched with the axiom $L_4$  given in Table \ref{lagrange} for $n=4$. 
\begin{table}[t]
\centering
\hrule
\begin{align*}
	(x+y)+z &= x + (y + z)\\
	x+y     &= y+x\\
	x+0     &= x\\
	x+(-x)  &= 0\\
	(x \cdot y) \cdot  z &= x \cdot  (y \cdot  z)\\
	x \cdot  y &= y \cdot  x\\
	1\cdot x &= x \\
	x \cdot  (y + z) &= x \cdot  y + x \cdot  z\\
	(x^{-1})^{-1} &= x \\
	x \cdot (x \cdot x^{-1}) &= x
\end{align*}
\hrule
\caption{The set $\Md$ of axioms for meadows}
\label{meadow}
\end{table}
\begin{table}[t]
\centering
\hrule
\begin{align*}
	(1+x_1^2 + \cdots + x_n^2)\cdot (1+x_1^2 + \cdots + x_n^2)^{-1}&=1\ \ \ (L_n)
\end{align*}
\hrule
\caption{The axiom schema $L_n$}
\label{lagrange}
\end{table}
The characterization above has been sharpened in \cite{BM2011} where it is shown that
\[
\mathbb{Q}_0 \cong \mathcal{I}(\Sigma_{\Md} , \Md+ L_2).\hfill{(\ddag)}
\]
 {\mbox{}  In \cite{BB15} it is proved that every finite specification of $\mathbb{Q}_0$ can be given in the form $\Md + e$|the meadow axioms enriched with a single equation. Moreover,  observe that both $L_4$ and $L_2$ do not hold in the presence of the imaginary unit $i$. In general, every finite algebraic specification of  $\mathbb{Q}_0$ contains an equation not valid in $\mathbb{C}_0$|the  zero-totalized expansion field of the complex numbers. Again, this is proved in \cite{BB15}.
 
In the sequel we denote by $(\mathbb{Z}/p\mathbb{Z})_0$ the zero-totalized expansion of the finite prime field of order $p$. Moreover, we define for $n\in \mathbb{N}$ the numerals $\underline{n}$ by $\underline{0}:= 0$ and $\underline{n+1}:=\underline{n}+1$. 
 A necessary and sufficient condition for  an initial algebra specification of $\mathbb{Q}_0$ is given in the following theorem.
\begin{theorem}\label{initialspec}
Let $E$ be a set of meadow equations. Then
\[
\mathbb{Q}_0 \cong \mathcal{I}(\Sigma_{\Md} , \Md+ E) \text{ if and only if  for all prime numbers $p$, } (\mathbb{Z}/p\mathbb{Z})_0 \not \models E.
\]
\end{theorem}
{\bf Proof}: Assume $\mathbb{Q}_0 \cong \mathcal{I}(\Sigma_{\Md} , \Md+ E)$. Then $\Md+E \vdash \underline{p}\cdot \underline{p}^{-1}=1$ for all prime numbers $p$. Hence $(\mathbb{Z}/p\mathbb{Z})_0 \not \models E$ for all primes $p$. For the converse, assume $(\mathbb{Z}/p\mathbb{Z})_0 \not \models E$ for all primes $p$. Recall that every  minimal meadow is a subdirect product of  minimal zero-totalized expansion fields (see e.g.\ \cite{BR10}). The minimal zero-totalized expansion fields are the zero-totalized prime fields  $(\mathbb{Z}/p\mathbb{Z})_0$ and $\mathbb{Q}_0$. In particular,  the initial algebra of $\Md+ E$ is such a subdirect product. Since every  $\mathbb{Z}/p\mathbb{Z}_0$ is not a model of $E$, it follows that $\mathcal{I}(\Sigma_{\Md} , \Md+E )$ must be isomorphic to $\mathbb{Q}_0$. \hfill $\Box$

An application of this theorem yielding a positive result is given below. First we  recall a few facts from the theory of numbers (see e.g.\ \cite{D52}, Ch.\ 3). For every odd prime $p$ one of the two congruences $p\equiv 1 \mod 4$ or $p\equiv 3 \mod 4$ hold.
Given a prime $p$ and a natural number $0<n<p$, $n$ is called a {\em quadratic residue} of $p$ if there exists a natural number $x$ such that $x^2\equiv n \mod p$. If the congruence is insoluble, $n$ is said to be a {\em quadratic non-residue}. Every prime $p$ has quadratic residues since $1^2 \equiv 1 \mod p$, but for $p>3$ there are more: e.g. if $p=19$, the quadratic residues are 1, 4, 5, 6, 7, 9, 11, 16 and 17. In general, every odd  prime $p$ has $\frac{p-1}{2}$ quadratic residues. If $p\equiv 1 \mod 4$ then the lists of quadratic residues and quadratic non-residues are both symmetrical in the sense that if $n$ is a (non-)quadratic residue then $p-n$ is one. On the other hand, if $p\equiv 3 \mod 4$, then $n$ is a quadratic residue if and only if $p-n$ is a quadratic non-residue (see e.g.\ the case of $p=19$). The quadratic residues and non-residues have a simple \emph{multiplicative property}: the product of two residues or of two non-residues is a residue.
\begin{example}
There also exist finite initial algebra specifications of $\mathbb{Q}_0$  of the form $\Md+ e$ where $e$ is a single variable equation.
Consider the equation $f(x)\cdot  f(x)^{-1}=1$ where $f(x)= (x^2-2)(x^2-3)(x^2-6)$. Inspection shows that $f(x)$ has no rational root. Thus $\mathbb{Q}_0\models f(x)\cdot  f(x)^{-1}=1$. On the other hand, $f(x)$ has a root modulo every prime number $p$: 
\begin{itemize}
\item If $p=2$ then $f(x)$ has a root  modulo $p$ for $x=0$.
\item If $p>2$ we apply the multiplicative property of quadratic residues and non-residues. If $(x^2 - 2)$ or $(x^2 -3)$ have a root modulo $p$, then $f(x)$ has a root modulo $p$. If neither $(x^2-2)$ nor $(x^2-3)$ has  a root modulo $p$ then both $2$ and $3$ are non-residues of $p$. Hence $6$ is a residue of $p$, i.e.\ $(x^2-6)$ has a root modulo $p$ and thus $f(x)$ has a root modulo $p$.
\end{itemize}
It follows that $ (\mathbb{Z}/p\mathbb{Z})_0 \not \models f(x)\cdot  f(x)^{-1}=1$ for every prime number $p$. By the above theorem we may conclude that $\Md + f(x)\cdot  f(x)^{-1}=1$ is an initial algebra specification of $\mathbb{Q}_0$.
\end{example}
 In the next section we apply Theorem \ref{initialspec} in order to give a negative result.
 
 \section{A negative result} 
 A general question concerning the specification of the rationals is whether there exists a logical weakest initial algebra specification. In this section we show that the weakening from $L_4$ to $L_2$ cannot be prolonged in a straightforward way.
 
 \begin{table}[h]
\centering
\hrule
\begin{align*}
	(1+\underline{n} + x^2)\cdot (1+\underline{n} + x^2)^{-1}&=1	\ \ \ (H_n)\\
         (1+\underline{n})\cdot (1+\underline{n})^{-1}&=1\ \ \ (C_n)
\end{align*}
\hrule
\caption{The axiom schemas $H_n$ and $C_n$}
\label{lagrange2}
\end{table}
Substituting 0 for $x$ in $H_n$  in Table \ref{lagrange2},  we obtain the axiom $C_n$. In \cite{BT2007} it is proved that 
$\{C_n\mid n\in \mathbb{N}\}$  together with $\Md$ specify the rational numbers.  The question then arises whether $\mathbb{Q}_0$ can be specified by $\Md + \Gamma$ for some finite subset  $\Gamma\subset \{H_n\mid n\in \mathbb{N}\}$. We give a negative answer below. 

Consider the function $f:\mathbb{N} \rightarrow \mathbb{N}$ defined by
\[
f(n)=
\begin{cases}
0 & \text{if $n\leq 1$ or $n$ is composite},\\
n-max\{i \mid \text{$i$ is a quadratic residue of $n$}\}& \text{if $n$ is prime}. 
\end{cases}
\]
E.g.\ $f(19)=2$.  In \cite{OEIS}, a table for the values of $f(n)$ can be found for the first $10^5$ natural numbers (see entry A088192). The  occurring  values increase very slowly: the largest value found in that table is 43.

In \cite{W07}, Wright proved the following theorem (see Theorem 2.3).
\begin{theorem}
Every nonempty finite subset of $\mathbb{N}^+$ is a set of quadratic residues for infinitely many primes.
\end{theorem}
As a corollary we obtain that the function $f$ is unbounded.
\begin{corollary}\label{unbounded}
$f$ is unbounded.
\end{corollary}
{\bf Proof}: For $n\in \mathbb{N}$ with $n>2$ consider the set $A=\{1,2,3, \ldots , n\}$. By the previous theorem we can pick a prime $p$ such that every element of $A$ is a quadratic residue of $p$. In particular, 2 is a quadratic residue of $p$. It follows from Gauss's lemma that $p\equiv 7 \mod 8$ and hence $p\equiv 3 \mod 4$. Thus, since $1, 2, 3, \ldots , n$ are all quadratic residues, $p-1, p-2, p-3, \ldots, p-n$ are all quadratic non-residues. So $max\{i \mid \text{$i$ is a quadratic residue of $p$}\}< p-n$. Therefore $f(p)>n$. \qed

We now prove that for every finite $\Gamma\subset \{H_n\mid n\in \mathbb{N}\}$ there exists a prime $p$ with $(\mathbb{Z}/p\mathbb{Z})_0\models \Gamma$.
\begin{proposition}
Let $\Gamma \subset  \{H_n\mid n\in \mathbb{N}\}$ be finite. Then there exists a prime $p$ such that $(\mathbb{Z}/p\mathbb{Z})_0\models \Gamma$.
\end{proposition}
{\bf Proof}: Say $\Gamma=\{H_0, \ldots , H_n\}$. Pick a prime $p$ such that $f(p)>n+1$. Suppose $(\mathbb{Z}/p\mathbb{Z})_0\not \models H_m$ for some $0\leq m\leq n$. We derive a contradiction as follows. By the assumption, there exists $0\leq x <p$ with $1+m +x^2=0$, i.e.\ $x^2 \equiv p-(m+1) \mod p$. Hence $p-(m+1)$ is a quadratic residue of $p$. Therefore 
$p-(m+1) \leq max\{i \mid \text{$i$ is a quadratic residue of $p$}\}$ and hence
\[
m+1 = p-(p-(m+1)) \geq p- max\{i \mid \text{$i$ is a quadratic residue of $p$}\}=f(p) > n+1\geq m+1. \hfill \qed
\]
We can now show that $\mathbb{Q}_0$ cannot be specified by a finite set of $H_n$-axioms.
\begin{theorem}
Let $\Gamma \subset  \{H_n\mid n\in \mathbb{N}\}$ be finite. Then $\mathbb{Q}_0 \not  \cong \mathcal{I}(\Sigma_{\Md} , \Md+ \Gamma)$.
\end{theorem}
{\bf Proof}:  Immediately by the preceding proposition and Theorem \ref{initialspec}. \hfill{\qed}

\begin{corollary}
$\mathbb{Q}_0 \not \cong \mathcal{I}(\Sigma_{\Md} , \Md+ L_1)$
\end{corollary}

{\bf Remark:} Observe that $\mathcal{I}(\Sigma_{\Md} , \Md+ L_1)$ is a non-cancellation meadow of characteristic 0 which does not contain $\mathbb{Q}_0$ as a subalgebra  since $\mathbb{Q}_0 \not \cong \mathcal{I}(\Sigma_{\Md} , \Md+ L_1)$. The existence of such a structure  has already been proved in \cite{BBP15} (see  Theorem 2.1). However, the proof given there relies on the compactness theorem.
\section{An open question}
Open questions arise when we extend the rational numbers.
E.g.\ consider the meadow of Gaussian rationals|denoted $\mathbb{Q}_0(i)$|obtained by adjoining the imaginary number $i$ to the meadow of rationals. It is not difficult to see that the polynomial given in the example in Section 1 also yields an initial algebra specification of the Gaussian rationals. Since  $f(x)=(x^2-2)(x^2 - 3)(x^2-6)$ has only real roots, there exist no complex rational ones. It follows  that $\mathbb{Q}_0(i)\models f(x)\cdot f(x)^{-1}=1$. Moreover, working in $\Md + \{i^2+1 =0\}$ it can be shown that every closed term over $\Sigma_\Md \cup \{i\}$ is provable equal to a term of the form $l\cdot m^{-1} +p\cdot q^{-1} \cdot i$ where $l,m,p,q$ are numerals. Hence $\Md + \{f(x)\cdot f(x)^{-1}=1, i^2+1 =0\}$ yields an initial algebra specification of $\mathbb{Q}_0(i)$. 

Up to isomorphism, there exists only one simple extension of the rational numbers generated by a transcendental number $t$.
We are unable to prove or disprove the existence of a finite initial algebra specification of $\mathbb{Q}_0(t)$ with $t$ a new constant interpreted as a trancendental element.
}

{\bf Acknowledgement}: We are indebted to Rob Tijdeman for valuable discussions including a proposal for proving Corollary \ref{unbounded} and the suggestion to make use of \cite{W07}.

\end{document}